\documentclass[preprint,12pt]{elsarticle}

\usepackage{graphics}

\usepackage{graphicx}

\usepackage{epsfig}
\usepackage{dsfont}

\usepackage{amssymb,amsthm}
\usepackage{pifont,natbib}

\usepackage[nodots,nocompress]{numcompress}

\journal{JCAM}
\begin{document}
\begin{frontmatter}


\title{Estimation of Hurst Parameter of Fractional Brownian Motion Using CMARS Method}


\author[bir]{F.˜Yerlikaya-\"Ozkurt\corref{cor1}}
\ead{fatmayerlikaya@gmail.com}

\author[iki]{C.˜Vardar-Acar}
\ead{cvardar@etu.edu.tr}

\author[bir]{Y.˜Yolcu-Okur}
\ead{yyolcu@metu.edu.tr}
\author[bir]{G.-W.˜Weber}
\ead{gweber@metu.edu.tr}

\cortext[cor1]{Corresponding author}

\address[bir]{Institute of Applied Mathematics, METU, 06800 Ankara, Turkey}
\address[iki]{Department of Mathematics, TOBB ETU, 06530 S\"o\c g\"ut\"oz\"u, Ankara, Turkey}

\begin{abstract}
%Stochastic Differential Equations (SDEs) generated by
%fractional Brownian motion (fBm) with Hurst parameter,
%$H$, are widely used to represent noisy and real-world problems. They
%play an important role in many fields of science such as finance,
%physics, biotechnology, meteorology and engineering. The crucial parameter in such a kind of SDEs is the Hurst parameter. The Hurst parameter $H$ of fBm explains the dependence structure of the data.

In this study, we develop a new theory of estimating Hurst parameter using conic multivariate adaptive regression splines (CMARS) method. We concentrate on the strong solution of stochastic differentional equations (SDEs) driven by fractional Brownian motion (fBm). The superiority of our approach to the others is, it not only estimates the Hurst parameter but also finds spline parameters of the stochastic process in an adaptive way. We examine the performance of our estimations using simulated test data.

\end{abstract}

\end{frontmatter}

\section{Introduction}
\textit{Fractional Brownian motion (fBm)} is a widely used concept for modelling various
situations such as the level of water in a river, the temperature at a specific place,
empirical volatility of a stock, the price dynamics of electricity. It
appears naturally in these phenomena because of its capability of
explaining the \textit{dependence} structure in real-life observations. This
structure in fBm is represented by its Hurst parameter
$H$. A fBm with Hurst parameter $H>1/2$ is called a \textit{persistent} process,
i.e., the increments of this process are positively correlated. On
the other hand, the increments of a fBm with $H<1/2$ is called an
\textit{anti-persistent} process with increments being negatively correlated.
For $H=1/2$, fBm corresponds to Brownian motion which has
independent increments. For further information on fBm and its applications, see \cite{oksendal, mishura, rao}.

It is highly important to identify the value of Hurst parameter in order to understand the structure of the process and its applications since the calculations dramatically differ according to the value of $H$. Therefore, some techniques have been developed to estimate Hurst parameter which can be categorized into three groups; heuristics, maximum likelihood and wavelet-based estimators. In the group of heuristics estimators, there is R/S estimator which was firstly proposed by Hurst \cite{hurst}, followed by the methods of correlogram, variogram, variance plot, and partial correlations plot. Due to lack of accuracy of heuristics estimators, maximum likelihood estimators (mle) were developed. Being weakly consistent is the main disadvantage of mle. In paralel to mle, wavelet-based estimators were suggested because of the popularity of wavelet decomposition of fBm \cite{viens, matos}.

In search of faster and efficient ways to estimate the Hurst parameter $H$,
we suggest a new numerical and computational way, \textit{conic multivariate adaptive regression splines (CMARS)}. CMARS is an alternative approach to the well-known data mining tool
\emph{multivariate adaptive regression splines (MARS)}. It is based
on a penalized residual sum of squares (PRSS) for MARS as a Tikhonov
regularization (TR) problem. CMARS treats this problem by a
continuous optimization technique, in particular, the framework of
\emph{conic quadratic programming (CQP)}. These convex optimization
problems are very well-structured, herewith resembling linear
programs and, hence, permitting the use of \textit{interior point methods}.

This paper is organized as follows; in Section 2, we start with explaning the properties of our madel given as, SDEs driven by fBm. In Section 3, we introduce the method CMARS relating it to the Hurst parameter estimation of our model. In Section 4, we give an application of our study, in order to test the theory we have developed. Finally, we present a brief conclusion and a general outlook of our study.

%%%%%%%%%%%%%%%%%%%%%%%%%%%%%%%%%%%%%%%%%%%%%%%%%%%%%%
\section{Stochastic Differential Equations with Fractional Brownian Motion}
Stochastic Differential Equations (SDEs) generated by
fBm are widely used to represent noisy and real-world problems. They
play an important role in many fields of science such as finance,
physics, biotechnology and engineering. In this section, we briefly recall some concepts on fBm and stochastic differential equations driven by fBm.

\subsection{Fractional Brownian Motion}
Let $H$ be a constant in the interval $(0,1)$. FBm $(W^{H}(t))_{{t \geq
0}}$ with Hurst parameter $H$, is a continuous and centered Gaussian
process with covariance function
\[
E[W^{H}(t) W^{H}(s)]=\frac{1}{2}(t^{2H}+s^{2H}-|t-s|^{2H}).
\]
We note that, for $H=1/2$, fBm corresponds to a standard Brownian motion which has independent increments.
For a \textit{standard fBm}, $W^{H}(t)$:
\begin{itemize}
\item $W^{H}(0)=0$ and $E[W^{H}(t)]=0$ for all $t\geq 0$.
    \item $W^{H}$ has homogenous increments, i.e., $W^{H}(t+s)-W^{H}(s)$ has the same law as $W^{H}(t),$ for all $s,t\geq 0$.
    \item $W^{H}$ is a Gaussian process and $E[(W^{H}(t))^2]=t^{2H} ~(t\geq 0)$, for all $H\in(0,1)$.
    \item $W^{H}$ has continuous trajectories.
\end{itemize}
%%%%%%%%%%%%%
The Hurst parameter $H$ of fBm explains the dependency of data. Indeed, the correlation between increments for $s, t \geq 0$ can be obtained by;
\[
\begin{array}{c}
\mathds{E} [(W^H (t+h)-W^H(t))(W^H (s+h)-W^H(s))]=  \\
\frac{h^{2 H}}{2} [ (n+1)^{2 H} + (n-1)^{2 H} - 2 n^{2 H} ].
\end{array}
\]

\begin{figure}[h]
   \centering
      \includegraphics[angle=0,width=90mm]{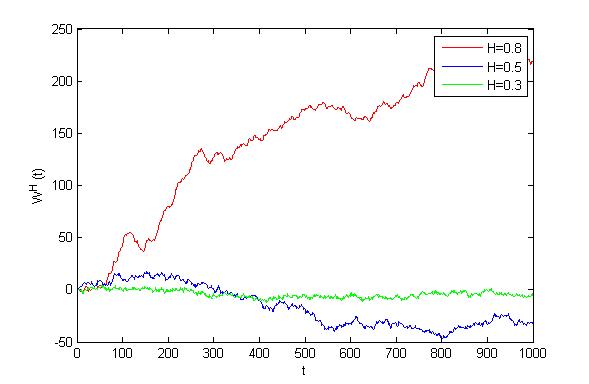}
       \caption{Sample paths of fBm with different Hurst parameter values simulated by Cholesky method.}
   \label{fig:fbm}
 \end{figure}

It can be seen that observations with $H>1/2$ have positively correlated increments and display long-range dependence, while the observations with
$H<1/2$ have a negatively correlated increments and display a short-range dependence structure (see Figure 1). Therefore, it is crucial to find the Hurst parameter of a stochastic process
for understanding the structural behaviour of this phenomena. In this study, we concentrate on finding $H$ for the stochastic processes which are the strong solutions of SDEs with fBm. Hence, we first
recall some fundamental properties of them.
%%%%%%%%%%%%%%%%%%%%%
\subsection{Stochastic Differential Equations Driven by Fractional Brownian Motion}
\par Suppose we have a stochastic process $ X= \{ X(t); ~ t \geq 0 \}$ defined on a filtered probability space $ (\Omega, \mathcal{F}, (\mathcal{F}_t)_{ t \geq 0}, \mathbb{P})$ which
is the strong solution of the following SDE:
\begin{equation} \label{sde1}
dX(t)=a (t, X (t))dt+ b(t, X(t))d W^H(t).
\end{equation}
Here, $ a(t, X(t))$ and $ b(t,X(t))$ are the drift and diffusion terms satisfying the conditions of existence and uniqueness theorem for $ t \geq 0$. Note that it is necessary to have the integrator as a semi-martingale in the theory of stochastic integration. However, since fBm is not a semi-martingale, one should extend the usual settings as in It\^o integral and define the integration with respect to fBm in a new pathwise integration technique. Alos et al. \cite{alos} construct the theory of integration with respect to general Gaussian proceses to overcome this. For further studies on this topic, see \cite{mishura, rao}.
%%%%%%%%%%%%%%%%%%%%%%%%%%%%
\par There have been comprehensive studies on statistical inferences for processes satisfying SDEs driven by Brownian motion. However, the recent interest is on SDEs driven by fBm since there have not been adequate studies on this topic. The purpose of this study is to estimate the Hurst parameter of the following SDE
\begin{equation} \label{sde2}
dX(t)=a(t,X(t))dt + b~dW^H (t), \qquad X_0 \in \mathbb{R},
\end{equation}
by \textit{Conic Multivariate Adaptive Regression Splines (CMARS)} methodology. Note that $b(t, X(t))\equiv b$ term in equation (\ref{sde1}) is taken as \textit{constant}.
%%%%%%%%%%%%%%%%%%%%%%%%%%%%%%%%%%%%%%%

\section{Estimation of Hurst Parameter Using Conic Multivariate Adaptive Regression Splines Method}
In this section, as an alternative to the existing methods of estimation of Hurst parameter, CMARS and the related proposed methodology will be introduced. For that purpose, firstly, we give a brief description of CMARS method and then we mention about the methodology and show how to apply this technique for finding the Hurst parameter of SDE defined in equation (\ref{sde2}).

\subsection{Method of Conic Multivariate Adaptive Regression Splines}
CMARS method is an alternative approach to the well-known data mining tool \textit{Multivariate Adaptive Regression Splines (MARS)}. It makes no specific assumption about the underlying functional relationship between the dependent and independent variables to estimate a general model function \cite{friedman1}. CMARS is introduced by linear combinations of the basis functions (BFs) that are used in MARS. The selection of BFs is data-based and specific to the problem at hand. CMARS uses one-dimensional BFs of the form $c^{+}(x,\tau)=[+(x-\tau)]_{+}$ and $c^{-}(x,\tau)=[-(x-\tau)]_{+}$, where $[q]_{+}:=\max \left\{0,q\right\}$ (see \cite{fatmatez, mainpaper} for further details). Each function is piecewise linear, with a knot at the value $\tau$, and the corresponding couple of function is called a \textit{reflected pair}. A set of BFs is given as follows:
\[
\wp:=\left\{(x_{j}-\tau)_{+}, (\tau-x_{j})_{+}\;|\;\tau \in \left\{ x_{1,j},x_{2,j},...,x_{N,j}\right\},\; j \in \left\{1,2,...,p \right\} \right\}.
\]

\begin{figure}[h]
    \centering
        \includegraphics[width=0.9\textwidth]{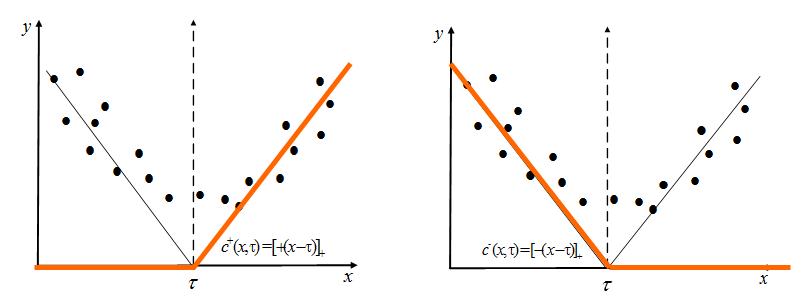}
    \caption{Visualization of the BFs.}
    \label{fig:BF}
\end{figure}

A CMARS model function $f$ is represented by a linear combination of BFs which is successively built up by the set $\wp$ as described below:
\begin{equation} \label{cmars}
Y=f(\textbf{\textit{x}})+ \epsilon=\theta_{0}+\sum^{M}_{m=1} \theta_{m} \psi_{m}(\textbf{\textit{x}}^{m})+ \epsilon.
\end{equation}

\noindent Here $Y$ is a response variable, $\textbf{\textit{x}}^{m}=(x_{1},x_{2},...,x_{p})^{T}$ a vector of predictors for the corresponding $m$th multivariate basis function. Furthermore, $\theta_{m}$ are the unknown coefficients for the $m$th basis function $(m=1,2,...,M)$ or for the constant 1 $(m=0)$, and  $\epsilon$ is an additive stochastic component which is assumed to have zero mean and finite variance. In equation (\ref{cmars}), $\psi_{m}$ $(m=1,2,...,M)$ are BFs as products of two or more one-dimentional BFs. Such interaction BFs are created by multiplying an existing basis function with a truncated linear function, involving a new variable. The form of the $m$th BF can be written as follows:

\begin{equation} \label{BF}
\psi_{m}(\textbf{\textit{x}}^{m}):=
{\prod^{K_{m}}_{j=1}}[s_{\kappa^{m}_{j}} \cdot (x_{\kappa^{m}_{j}}-\tau_{\kappa^{m}_{j}})]_{+}.
\end{equation}

\noindent Here, $\textbf{\textit{x}}^{m}$ is the vector of variable contributed to the $m$th
BF, $K_{m}$ is the number of truncated linear functions multiplied in the $m$th
BF, $x_{\kappa^{m}_{j}}$ is the predictor variable corresponding to the $j$th truncated
linear function in the $m$th BF, $\tau_{\kappa^{m}_{j}}$ is the knot
value corresponding to the variable $x_{\kappa^{m}_{j}}$,  and $s_{\kappa^{m}_{j}}$ is the selected
sign $+$1 or $-$1.

CMARS is constructed by a \textit{Penalized Residual Sum of Squares (PRSS)} parameter estimation problem, instead of an ordinary least-squares estimation problem as it occurs in MARS method. The PRSS problem aims at accuracy and a smallest possible complexity of the model. PRSS with penalty parameters $\lambda_{m}$ and with $M_{max}$ BFs which are accumulated in the first part of the MARS algorithm, has the following form:

\begin{equation} \label{prss}
PRSS:={{\sum^{N}_{i=1}}}(y_{i}-f(\tilde{\textbf{\textit{x}}}_{i}))^{2}+%
{{\sum^{M_{max}}_{m=1}}}\lambda_{m}%
{\sum^{2}_{\stackrel{\left| \mbox{\boldmath{\tiny{${\alpha}$}}}\right|=1}{\mbox{\boldmath{\tiny{${\alpha}$}}} =(\alpha_{1},\alpha_{2})^{T}}}}\ %
{\sum_{\stackrel{r < s}{r,s\in V_{m}}}} \int_{Q^{m}} \theta^{2}_{m}\left[D^{\mbox{\boldmath{\tiny{${\alpha}$}}}}_{r,s}\psi_{m}(\textbf{\textit{t}}^{m})\right]^{2}d\textbf{\textit{t}}^{m},
\end{equation}

\noindent where $V_{m}:=\left\{\kappa^{m}_{j}|j=1,2,...,K_{m}\right\}$ is the variable
set associated with the ${\it m}$th BF, $\textbf{\textit{t}}^{m}=\left(t_{m_{_{1}}},t_{m_{_{2}}},...,t_{m_{_{K_{m}}}}\right)^{T}$ represents
the vector of variables which contribute to the ${\it m}$th BF. Moreover, $D^{\mbox{\boldmath{\tiny{${\alpha}$}}}}_{r,s} \psi_{m}(\textbf{\textit{t}}^{m}):=
\frac{\partial^{\left|\mbox{\boldmath{\tiny{${\alpha}$}}}\right|} \psi_{m}}{\partial^{\alpha_{1}}t^{m}_{r}\;\partial^{\alpha_{2}}t^{m}_{s}}(\textbf{\textit{t}}^{m})$ for
${\mbox{\boldmath${\alpha}$}}=(\alpha_{1},\alpha_{2})^{T}$, $\left|{\mbox{\boldmath${\alpha}$}}\right|:=\alpha_{1}+\alpha_{2}$, where $\alpha_{1},\alpha_{2}\in\left\{0,1\right\}$.

After using the same penalty parameter $\lambda=\lambda_{m}$ for each derivative, PRSS turns into a \textit{Tikhonov regularization problem} as described below:

\begin{equation} \label{TR}
PRSS \approx \left\|\textbf{\textit{y}}-\mbox{\boldmath${\psi}$}(\tilde{\textbf{\textit{d}}})\mbox{\boldmath${\theta}$}\right\|^{2}_{2}+
\lambda \left\|\textbf{\textit{L}}\mbox{\boldmath${\theta}$}\right\|^{2}_{2},%
\end{equation}

\noindent where $\mbox{\boldmath${L}$}$ is constructed by the discretizations of the high-dimentional integrals given in equation (\ref{prss}). The model approximations as presented in equation (\ref{TR}) are carefully prepared. They play an important role in order to raise a final model approximation which is linear in the unknown spline parameters. After unifying some discretised complexity terms and including them into inequality constraints,
a \textit{Conic Quadratic Programming (CQP)} problem is obtained which uses \textit{interior point methods} \cite{nemirovski1, nemirovski2}. The formulation of CQP is given as follows:

\begin{eqnarray} \label{CQP}
\min_{t,\mbox{\boldmath{\tiny{${\theta}$}}}} \ \ \ \ \ t, \quad\quad \quad\quad \quad  \quad \quad\quad \quad\quad \quad  \quad  \quad\quad \quad\quad \quad \nonumber\\
\textnormal{subject \ to} \ \ \left\|\textbf{\textit{y}}-\mbox{\boldmath${\psi}$}(\tilde{\textbf{\textit{d}}})\mbox{\boldmath${\theta}$}\right\|_{2} \leq t, \quad \left\|\textbf{\textit{L}}\mbox{\boldmath${\theta}$}\right\|_{2} \leq \sqrt{\widetilde{M}},
\end{eqnarray}

\noindent referring to some chosen complexity bound $\widetilde{M}\geq 0$.

\subsection{Discretization of Stochastic Differential Equations with Fractional Brownian Motion}
In general, the distribution of the stochastic process $ \left\{X(t); ~ t \geq 0 \right\}$ is not known. Therefore, the discretized version of the SDEs, $(\hat{X}_{i})_{i \in \mathbb{N}}$, should be simulated \cite{kloeden}. There are many discretization schemes for the SDEs generated by fBm such as  Euler and Milstein Scheme (see \cite{milstein} for further details). In this study, \textit{Euler approximation} is used since Milstein approximation contains the derivatives of the diffusion term in equation (\ref{sde2}) which is equal to zero. The Euler approximation of the equation (\ref{sde2}) is:

\begin{eqnarray}
\hat{X}_{i+1}=\hat{X}_{i}+a(\hat{X}_{i},t_{i})(t_{i+1}-t_{i})+b(\hat{X}_{i},t_{i})(W^{H}_{i+1}-W^{H}_{i}) \quad (i \in \mathbb{N}).
\end{eqnarray}

\noindent For finitely many given data points $(\bar{X}_{i},\bar{t}_{i})$ $(i=1,2,\ldots,N)$, the symbolic form of the \textit{approximation} can be given as follows:

\begin{eqnarray} \label{euler}
\dot{\bar{X}}_{i}=a(\bar{X}_{i},\bar{t}_{i})+b(\bar{X}_{i},\bar{t}_{i})\frac{\Delta W^{H}_{i}}{\bar{h}_{i}},
\end{eqnarray}

\noindent where $\Delta W^{H}_{i}=W^{H}_{i+1}-W^{H}_{i}$ is a centered Gaussian random variable, $\bar{h}_{i}=\bar{t}_{i+1}-\bar{t}_{i}:=\Delta\bar{t}_{i}$ represents step lengths and
\[
\dot{\bar{X}}_{i}=\left\{
\begin{array}{lr}
\frac{\bar{X}_{i+1}-\bar{X}_{i}}{\bar{h}_{i}} \ \ \ \ \ \ \textnormal{if} \ \ \ \ \ \ i=1,2,...,N-1, \\
\\
\frac{\bar{X}_{N}-\bar{X}_{N-1}}{\bar{h}_{N}} \ \ \ \ \textnormal{if} \ \ \ \ \ \ i=N,
\end{array}
\right.
\]

\noindent represents the difference quotients raised on the $i$th data value. A more compact form of the equation (\ref{euler}) is defined by
\begin{eqnarray} \label{abb}
\dot{\bar{X}}_{i}=\bar{G}_{i}+\bar{F}_{i}c_{i},
\end{eqnarray}

\noindent where $\bar{G}_{i}:=a(\bar{X}_{i},\bar{t}_{i})$, $\bar{F}_{i}:=b(\bar{X}_{i},\bar{t}_{i})$, and $c_{i}:= \Delta W^{H}_{i}/\bar{h}_{i}$. Note that equation (\ref{abb}) can be considered as an approximation of the problem. The expressions stated until the end of Section 3 are described parametrically with respect to the Hurst parameter $H$. In Section 4, we shall specify it by numeric values.

%%%%%%%%%%%%%%%%%%%%%%%%%%%%%%%%%%
\subsection{Parameter Estimation}

To determine the unknown values in equation (\ref{abb}), the following minimization problem is constructed using some abbreviated notation of the approximation \cite{pakize}:

\begin{eqnarray}
\min_{\mbox{\boldmath{\tiny{${\theta}$}}}} \ \ \ \ \ {{\sum^{N}_{i=1}}} \left\|\dot{\bar{X}}_{i}-\left(\bar{G}_{i}+\bar{F}_{i}c_{i}\right)\right\|^{2}_{2}. \nonumber
\end{eqnarray}

\noindent Here, $\mbox{\boldmath{${\theta}$}}$ comprises all unknown parameters in the Euler approximation. To solve this optimization problem and to give a smoother, regularized  approximation to the data, we employ CMARS method which controls any high ``variation'' in the data. CMARS' BFs are gradually constructed for the approximation of $\bar{G}_{i}$ and $\bar{F}_{i}$ with data $\bar{\textbf{\textit{U}}}^{l}_{i,B}, ~ \bar{\textbf{\textit{U}}}^{m}_{i,C}=(\bar{X}_{i},\bar{t}_{i})$ according to the following approaches \cite{fatmaSDE}:

\[
\bar{G}_{i}= \alpha_{0}+{{\sum^{d^{B}}_{l=1}}} \alpha_{l} B_{l}(\bar{\textbf{\textit{U}}}^{l}_{i,B}), \quad \textnormal{and} \quad \bar{F}_{i}c_{i}=\beta_{0}+{{\sum^{d^{C}}_{m=1}}} \beta_{m} C_{m}(\bar{\textbf{\textit{U}}}^{m}_{i,C}).
\]

\noindent Here, the forms of the BFs are $B_{l}(\bar{\textbf{\textit{U}}}_{B}^{l})={\prod^{2}_{k=1}}[s_{\kappa^{l}_{k}}^{B} \cdot (x_{\kappa^{l}_{k}}^{B}-\tau_{\kappa^{l}_{k}}^{B})]_{+}$ and $C_{m}(\bar{\textbf{\textit{U}}}_{C}^{m})={\prod^{2}_{k=1}}[s_{\kappa^{m}_{k}}^{C} \cdot (x_{\kappa^{m}_{k}}^{C}-\tau_{\kappa^{m}_{k}}^{C})]_{+}$ as we described in equation (\ref{BF}). Here, we choose the numbers $K^{B}_{l}$ and $K^{C}_{m}$ (in sense of equation (\ref{BF})) as maximal, namely, as 2.

We construct the penalized residual sum of squares (PRSS) for our minimization problem in the following form:

\begin{eqnarray} \label{prss2}
PRSS&:=&{{\sum^{N}_{i=1}}}(\dot{\bar{X}}_{i}-\left(\bar{G}_{i}+\bar{F}_{i}c_{i}\right))^{2}+ \nonumber \\
&\ & {{\sum^{d^{B}}_{l=1}}}\lambda_{l}%
{\sum^{2}_{\stackrel{\left| \mbox{\boldmath{\tiny{${\alpha}$}}}\right|=1}{\mbox{\boldmath{\tiny{${\alpha}$}}} =(\alpha_{1},\alpha_{2})^{T}}}}\ %
{\sum_{\stackrel{r < s}{r,s\in V_{l}}}} {\int_{Q^{B}_{l}}} \alpha^{2}_{l}\left[D^{\mbox{\boldmath{\tiny{${\alpha}$}}}}_{r,s}B_{l}(\bar{\textbf{\textit{U}}}^{l}_{B})\right]^{2}d\textbf{\textit{U}}^{l}_{B}+ \nonumber \\
&\ & {{\sum^{d^{C}}_{m=1}}}\mu_{m}%
{\sum^{2}_{\stackrel{\left| \mbox{\boldmath{\tiny{${\alpha}$}}}\right|=1}{\mbox{\boldmath{\tiny{${\alpha}$}}} =(\alpha_{1},\alpha_{2})^{T}}}}\ %
{\sum_{\stackrel{r < s}{r,s\in V_{m}}}} {\int_{Q^{C}_{m}}} \beta^{2}_{m}\left[D^{\mbox{\boldmath{\tiny{${\alpha}$}}}}_{r,s}C_{m}(\bar{\textbf{\textit{U}}}^{m}_{C})\right]^{2}d\textbf{\textit{U}}^{m}_{C}.
\end{eqnarray}

\noindent Here, the multipliers $\lambda_{l}, ~ \mu_{m} \geq 0$ are smoothing parameters and they provide a tradeoff between both accuracy and complexity. To approximate two multi-dimensional integrals in equation (\ref{prss2}), parallepipes $Q^{B}_{l}=\left[a^{l}_{1,B},b^{l}_{1,B}\right]\times\left[a^{l}_{2,B},b^{l}_{2,B}\right]={\prod^{2}_{k=1}}Q^{l}_{k,B}$ and $Q^{C}_{m}=\left[a^{m}_{1,C},b^{m}_{1,C}\right]\times\left[a^{m}_{2,C},b^{m}_{2,C}\right]={\prod^{2}_{k=1}}Q^{m}_{k,C}$ which encompass all our input data are constructed. Then, the following discretization is applied for the first multi-dimensional integral:

\begin{eqnarray}
{\int_{Q^{B}_{l}}} \alpha^{2}_{l}\left[D^{\mbox{\boldmath{\tiny{${\alpha}$}}}}_{r,s}B_{l}(\bar{\textbf{\textit{U}}}^{l}_{B})\right]^{2}d\textbf{\textit{U}}^{l}_{B} &\approx & {{\sum^{(N+1)^{2}}_{i=1}}} \left(
{\sum^{2}_{\stackrel{\left| {\mbox{\boldmath{\tiny{${\alpha}$}}}}\right|=1}{{\mbox{\boldmath{\tiny{${\alpha}$}}}}  =(\alpha_{1},\alpha_{2})^{T}}}}\ %
{\sum_{\stackrel{r < s}{r,s\in V_{m}}}} \alpha^{2}_{l}
\left[D^{\mbox{\boldmath{\tiny{${\alpha}$}}}}_{r,s}B_{l}(\hat{\textbf{\textit{U}}}^{l}_{i,B})\right]^{2}\right)\Delta \hat{\textbf{\textit{U}}}^{l}_{i,B} \nonumber \\
&=& {{\sum^{(N+1)^{2}}_{i=1}}} \left(\bar{L}^{B}_{il}\right)^{2}\alpha^{2}_{l} =\left\|\bar{\textbf{\textit{L}}}^{B}_{l}\alpha_{l}\right\|^{2}_{2}.
\end{eqnarray}

\noindent The same discretization is also applied for the second multi-dimensional integral in equation (\ref{prss2}). For simplicity, we introduce PRSS in the following matrix notation:

\begin{eqnarray}
PRSS \approx \left\|\dot{\bar{\textbf{X}}}-\bar{\textbf{\textit{A}}}\mbox{\boldmath${\theta}$}\right\|^{2}_{2}+ {{\sum^{d^{B}}_{l=1}}}\lambda_{l}\left\|\bar{\textbf{\textit{L}}}^{B}_{l}\alpha_{l}\right\|^{2}_{2}+{{\sum^{d^{C}}_{m=1}}}\mu_{m}\left\|\bar{\textbf{\textit{L}}}^{C}_{m}\beta_{m}\right\|^{2}_{2},
\end{eqnarray}

\noindent where $\dot{\bar{\textbf{X}}}=\left(\dot{\bar{X}}_{1},\dot{\bar{X}}_{2},...,\dot{\bar{X}}_{N}\right)^{T}$,  $\bar{\textbf{\textit{A}}}=\left(\bar{\textbf{\textit{A}}}^{T}_{1},\bar{\textbf{\textit{A}}}^{T}_{2},...,\bar{\textbf{\textit{A}}}^{T}_{N}\right)^{T}$, $\mbox{\boldmath${\theta}$}=\left(\mbox{\boldmath${\alpha}$}^{T},\mbox{\boldmath${\beta}$}^{T}\right)^{T}$, $\mbox{\boldmath${\alpha}$}=\left(\alpha_{0},\alpha_{1},\alpha_{2},..., \alpha_{d^{B}}\right)^{T}$ and $\mbox{\boldmath${\beta}$}=\left(\beta_{0},\beta_{1},\beta_{2},..., \beta_{d^{C}}\right)^{T}$, $\bar{G}_{i}+\bar{F}_{i}c_{i}=\bar{\textbf{\textit{A}}}_{i}\mbox{\boldmath${\theta}$}$, $\bar{\textbf{\textit{L}}}^{B}_{l}=\left(L^{B}_{1l},L^{B}_{2l},...,L^{B}_{(N+1)^{2}l}\right)^{T}$ and $\bar{\textbf{\textit{L}}}^{C}_{m}=\left(L^{C}_{1m},L^{C}_{2m},...,L^{C}_{(N+1)^{2}m}\right)^{T}$ for  $l=1,2,...,d^{B}$ and  for $m=1,2,...,d^{C}$, respectively \cite{fatmaSDE}.

Using uniform penalization by taking the same $\lambda$ for each derivative term, the regularized approximation problem of PRSS turns into a \textit{Tikhonov regularization (TR)} problem:

\begin{eqnarray}
PRSS \approx \left\|\dot{\bar{\textbf{X}}}-\bar{\textbf{\textit{A}}}\mbox{\boldmath${\theta}$}\right\|^{2}_{2}+ \lambda \left\|\bar{\textbf{\textit{L}}}\mbox{\boldmath${\theta}$}\right\|^{2}_{2}.
\end{eqnarray}

\noindent Here, $\lambda=\lambda_{1}=...=\lambda_{d^{B}}=\mu_{1}=...=\mu_{d^{C}},$ and $\bar{\textbf{\textit{L}}}$ is an $(M_{max}+1)\times(M_{max}+1)$-diagonal matrix with first column $\textbf{\textit{L}}_{0}=\textbf{0}_{(N+1)^{2}}$ and the
other columns being the vectors $\textbf{\textit{L}}^{B}_{l}$, $\textbf{\textit{L}}^{C}_{m}$,   introduced above, where $M_{max}=d^{B}+d^{C}$.

As we just mentioned in Subsection 3.1, TR problem can be solved by a CQP program as given in equation (\ref{CQP}). In order to write the optimality condition for this problem, we firstly reformulate our program  as the subsequent \textit{primal problem}:

\begin{eqnarray} \label{CQP2}
&\ &\min_{t,\boldmath{\mbox{{\tiny${\theta}$}}}}  \ \ t,  \nonumber\\
&\ & \textnormal{such \ that } \ \
\mbox{\boldmath$\chi$} :=
\left[\begin{array}{cc}
\textbf{0}_{N} & \bar{\textbf{\textit{A}}} \\
\\
1 & \textbf{0}^{T}_{M_{max}+1} \end{array}\right]
\left[\begin{array}{c}  t \\
\\
\mbox{\boldmath$\theta$} \end{array}\right]+
\left[\begin{array}{c}   -\dot{\bar{\textbf{X}}} \\
\\
  0 \end{array} \right], \nonumber\\
&\ & \mbox{\boldmath$\eta$} := \left[ \begin{array}{cc}
\textbf{0}_{M_{max}+1} & \bar{\textbf{\textit{L}}} \\
\\
0 & \textbf{0}^{T}_{M_{max}+1} \end{array} \right] \left[\begin{array}{c} t \\
\\
 \mbox{\boldmath$\theta$} \end{array} \right]+
\left[\begin{array}{c}  \textbf{0}_{M_{max}+1} \\
\\
\sqrt{\widetilde{M}} \end{array} \right], \nonumber\\
&\ & \mbox{\boldmath$\chi$} \in \textit{L}^{N+1},\; \mbox{\boldmath$\eta$} \in \textit{L}^{M_{max}+2}.
\end{eqnarray}

\noindent Here, $\textit{L}^{N+1}$, $\textit{L}^{M_{max}+2}$ are the $(N+1)$- and $(M_{max}+2)$-dimensional \textit{ice-cream} (or \textit{second-order}) \textit{cones} \cite{nemirovski1}. The \textit{dual problem} to the latter problem is given by \\
\\

\begin{eqnarray}
&\ & \max \ \ \ \ \ (\dot{\bar{\textbf{X}}}^{T},0)\mbox{\boldmath${\omega}$}_{1}+\left(\textbf{0}^{T}_{M_{max}+1},- \sqrt{\widetilde{M}}\right)\mbox{\boldmath${\omega}$}_{2}, \ \ \ \ \ \ \nonumber\\
&\ & \left[\begin{array}{cc}
\textbf{0}^{T}_{N} & 1 \\
\\

\bar{\textbf{\textit{A}}}^{T} & \textbf{0}^{T}_{M_{max}+1}
 \end{array} \right] \mbox{\boldmath${\omega}$}_{1}+ \left[\begin{array}{cc}
\textbf{\textit{0}}^{T}_{M_{max}+1} & 0 \\
\\
\bar{\textbf{\textit{L}}}^{T} & \textbf{0}_{M_{max}+1}
 \end{array} \right] \mbox{\boldmath${\omega}$}_{2} = \left[\begin{array}{c}
1 \\
\\
\textbf{0}_{M_{max}+1}
 \end{array} \right], \ \ \ \ \ \  \nonumber\\
&\ & \mbox{\boldmath${\omega}$}_{1} \in \textit{L}^{N+1},\; \mbox{\boldmath${\omega}$}_{2} \in \textit{L}^{M_{max}+2}. \ \ \ \ \ \
\end{eqnarray}

\noindent A \textit{primal-dual optimal solution} $(t,\mbox{\boldmath${\theta}$}, \mbox{\boldmath${\chi}$}, \mbox{\boldmath${\eta}$}, \mbox{\boldmath${\omega}$}_{1}, \mbox{\boldmath${\omega}$}_{2})$ is obtained when the \textit{optimality conditions} given in equation (\ref{OC}) are satisfied:

\begin{eqnarray} \label{OC}
&\ & \mbox{\boldmath${\chi}$} :=
\left[\begin{array}{cc}
\textbf{0}_{N} & \bar{\textbf{\textit{A}}} \\
\\

1 & \textbf{\textit{0}}^{T}_{M_{max}+1}%
\end{array} \right]
\left[\begin{array}{c}
t \\
\\
\mbox{\boldmath${\theta}$}
\end{array} \right]
+
\left[\begin{array}{c}
 -\dot{\bar{\textbf{X}}} \\
\\
0%
\end{array} \right], \nonumber\\
&\ & \mbox{\boldmath${\eta}$} :=
\left[ \begin{array}{cc}
\textbf{0}_{M_{max}+1} & \bar{\textbf{\textit{L}}} \\
\\
0 & \textbf{0}^{T}_{M_{max}+1}
\end{array} \right]
\left[ \begin{array}{c}
t \\
\\
\mbox{\boldmath${\theta}$}
\end{array} \right]
+
\left[ \begin{array}{c}
\textbf{0}_{M_{max}+1} \\
\\
\sqrt{\widetilde{M}}
\end{array} \right], \nonumber\\
&\ & \left[ \begin{array}{cc}
\textbf{0}^{T}_{N} & 1 \\
\\
\bar{\textbf{\textit{A}}}^{T} & \textbf{0}^{T}_{M_{max}+1}
\end{array} \right]
\mbox{\boldmath${\omega}$}_{1}+
\left[ \begin{array}{cc}
\textbf{0}^{T}_{M_{max}+1} & 0 \\
\\
\bar{\textbf{\textit{L}}}^{T} & \textbf{0}_{M_{max}+1}
\end{array} \right]
\mbox{\boldmath${\omega}$}_{2}=
\left[ \begin{array}{c}
1 \\
\\
\textbf{0}_{M_{max}+1}
\end{array} \right], \nonumber\\
&\ & \mbox{\boldmath${\omega}$}_{1}^{T} \mbox{\boldmath${\chi}$} =0,\ \; \mbox{\boldmath${\omega}$}_{2}^{T} \mbox{\boldmath${\eta}$} =0, \nonumber\\
&\ & \mbox{\boldmath${\omega}$}_{1} \in \textit{L}^{N+1} ,\;\  \mbox{\boldmath${\omega}$}_{2} \in \textit{L}^{M_{max}+2} , \;\
\mbox{\boldmath${\chi}$} \in \textit{L}^{N+1} ,\;\  \mbox{\boldmath${\eta}$} \in \textit{L}^{M_{max}+2} .
\end{eqnarray}

\section{Application and Results}
In order to test the theory developed in the previous section, we start with simulating stochastic process for a fixed Hurst parameter $H$ using Cholesky method \cite{dieker}. Now, our aim is to estimate the exact value of this Hurst parameter of the simulated data. We generate various stochastic processes which are the strong solution of SDEs driven by fBm with different Hurst parameters. Next, we construct CMARS model for each generated process to find the best fit. For the implementation of CMARS algorithm, BFs are built using Salford MARS$^{\tiny{\textregistered}}$ software program \cite{salford1} as in \cite{matlab, mainpaper}. The optimization problem given in equation (\ref{CQP2}) is solved by using interior point methods (IPMs) via the optimization software  MOSEK ~\cite{mosek, nemirovski2}. Finally, we examine the performances of CMARS fits according to well-known performance measures such as mean absolute error (MAE), mean squared error (MSE), correlation coefficient ($r$), multiple coefficient of determination ($R^{2}$), adjusted $R^{2}$ (Adj-$R^{2}$), and proportion of residuals within three sigma (PWI). The steps described above are applied for $H$=0.2, $H$=0.3, $H$=0.7, and $H$=0.8. The results of the applications are summarized in Table 1.

\begin{table}[h]%
\scriptsize{
\caption{CMARS performances for fBm generated by $H$=0.2, $H$=0.3, $H$=0.7, and $H$=0.8.}
\label{tab:2}       % Give a unique label
\begin{center}
\begin{tabular}
{|l|l|l|l|l|l|l|}\hline
&\multicolumn{6}{|c|}{\tiny{\textbf{Performance Measures}}} \\ \cline{2-7}
\tiny{\textbf{Hurst index}}  &  \tiny{\textbf{MAE}} & \tiny{\textbf{MSE}} & \tiny{\textbf{r}} & $\tiny{\textbf{R}^{2}}$ & $\tiny{\textbf {Adj-}\textbf{R}^{2}}$ & \tiny{\textbf{PWI}} \\ \hline

\tiny{$H=0.1$} & \tiny{0,8766} & \tiny{1,4827} & \tiny{0,1923} & \tiny{0,0370} & \tiny{-0,1651} & \tiny{1} \\\cline{2-7}

\tiny{\mbox{\boldmath${H=0.2}$}} & \tiny{0,6207*} & \tiny{0,7480*} & \tiny{0,9868*} & \tiny{0,9739*} & \tiny{0,9684*} & \tiny{1} \\\cline{2-7}

\tiny{$H=0.3$} & \tiny{0,8861} & \tiny{1,5266} & \tiny{0,0991} & \tiny{0,0098} & \tiny{-0,1979} & \tiny{1} \\\cline{2-7}

\tiny{$H=0.4$} & \tiny{0,8770} & \tiny{1,4733} & \tiny{0,2075} & \tiny{0,0430} & \tiny{-0,1577} & \tiny{1} \\ \cline{2-7}

\tiny{$H=0.5$} & \tiny{0,8839} & \tiny{1,5201} & \tiny{0,1281} & \tiny{0,0164} & \tiny{-0,1900} & \tiny{1} \\  \hline \hline

%\multicolumn{7}{|c|}{\small{\textbf{fBm generated by $\small{\textbf{H=0.3}}$}}} \\ \cline{1-7}
\tiny{$H=0.1$} & \tiny{0,7138} & \tiny{0,9776} & \tiny{0,2607} & \tiny{0,0679} & \tiny{-0,1276} & \tiny{1} \\\cline{2-7}
\tiny{$H=0.2$} &  \tiny{0,7162} & \tiny{0,9901} & \tiny{0,2400} & \tiny{0,0576} & \tiny{-0,1401} & \tiny{1} \\\cline{2-7}
\tiny{\mbox{\boldmath${H=0.3}$}} &  \tiny{0,3606*} & \tiny{0,2516*} & \tiny{0,984*} & \tiny{0,9699*} & \tiny{0,9636*} & \tiny{1} \\\cline{2-7}
\tiny{$H=0.4$} &  \tiny{0,6926} & \tiny{0,9387} & \tiny{0,3249} & \tiny{0,1055} &
\tiny{-0,0821} & \tiny{1} \\\cline{2-7}
\tiny{$H=0.5$} &  \tiny{0,7053} & \tiny{0,9763} & \tiny{0,2641} & \tiny{0,0697} & \tiny{-0,1254} & \tiny{1} \\\hline \hline

%\multicolumn{7}{|c|}{\small{\textbf{fBm generated by $\small{\textbf{H=0.7}}$}}} \\ \cline{1-7}
\tiny{$H=0.5$} &  \tiny{0,7031} & \tiny{0,9719} & \tiny{0,4743} & \tiny{0,2250} & \tiny{0,0623} & \tiny{1} \\\cline{2-7}
\tiny{$H=0.6$} &  \tiny{0,7048} & \tiny{0,9784} & \tiny{0,4688} & \tiny{0,2198} & \tiny{0,0560} & \tiny{0,9898} \\\cline{2-7}
\tiny{\mbox{\boldmath${H=0.7}$}} &  \tiny{0,3634*} & \tiny{0,2602*} & \tiny{0,9582*} & \tiny{0,9182*} & \tiny{0,9010*} & \tiny{1} \\\cline{2-7}
\tiny{$H=0.8$} &  \tiny{0,7041} & \tiny{0,9506} & \tiny{0,4919} & \tiny{0,2419} & \tiny{0,0828} & \tiny{1} \\\cline{2-7}
\tiny{$H=0.9$} &  \tiny{0,7081} & \tiny{0,9781} & \tiny{0,4691} & \tiny{0,2200} & \tiny{0,0563} & \tiny{1} \\\hline \hline

%\multicolumn{7}{|c|}{\tiny{\textbf{fBm generated by $\small{\textbf{H=0.8}}$}}} \\ \cline{1-7}
\tiny{$H=0.5$} &  \tiny{0,6068} & \tiny{0,7841} & \tiny{0,5914} & \tiny{0,3498} & \tiny{0,2133} & \tiny{1} \\\cline{2-7}
\tiny{$H=0.6$} &  \tiny{0,6359} & \tiny{0,8015} & \tiny{0,5783} & \tiny{0,3345} & \tiny{0,1948} & \tiny{1} \\\cline{2-7}
\tiny{$H=0.7$} &  \tiny{0,6053} & \tiny{0,7389} & \tiny{0,6176} & \tiny{0,3815} & \tiny{0,2517} & \tiny{1} \\\cline{2-7}
\tiny{\mbox{\boldmath${H=0.8}$}} &  \tiny{0,2009*} & \tiny{0,0822*} & \tiny{0,9883*} & \tiny{0,9768*} & \tiny{0,9720*} & \tiny{1} \\\cline{2-7}
\tiny{$H=0.9$} &  \tiny{0,6006} & \tiny{0,7294} & \tiny{0,6240} & \tiny{0,3894} & \tiny{0,2613} & \tiny{1} \\\hline
\end{tabular}
\end{center}
\tiny{* indicates better performance}}
\end{table}
In the case for anti-persistent processes, namely, $H=0.2$, $H=0.3$, the values of MAE and MSE are lower and the values of $R^{2}$, Adj-$R^{2}$, PWI and $r$ are higher than the values for the other Hurst parameter values. Similar results are also obtained for the case of persistent processes. Hence, this shows that according to performance measures criteria, the best CMARS fit gives us the correct Hurst parameter value.

\section{Conclusion and Outlook}

Recent developments in computer science provide environments in order to collect numerous data from various sources. Data mining methods enable us to analyze data for different purposes in many fields, such as finance, environment, and energy. One of the modern method of data mining, CMARS, has been developed as an alternative to the backward stepwise part of the MARS algorithm (see \cite{mainpaper}).

This paper gave a new contribution to Hurst parameter estimation theory for the strong solution of SDEs driven by fBm using CMARS technique. The main superiority of our approach to the others is that it not only estimates the Hurst parameter but it also finds spline parameters of the stochastic process. What is more, our representation of financial and other processes is empowered by all the modeling and numerical advantages of CMARS. By this, a bridge has been offered between convex optimization and Hurst parameter estimation theory.

In this pioneering paper, we followed a {\it two-level} approach with the determination of the parameters at the {\it lower level}, except of the Hurst-parameter which was chosen at the following {\it upper level}. This approach can be regarded as a {\it parametric optimization} (cf. \cite{guddat, jongen}). In future research, we will deepen and extend this approach by both more {\it model-free} strategies (e.g., from statistics and data mining), especially, more {\it model-based} ones, and with a comparison of them. The model-based approaches will be of a more {\it integrated} mathematical nature and in the analytical line that we initiated in this work. Through these investigations we intend to further contribute to a deeper understanding of our modern financial markets and to offer helpful mathematical decision tools for them. \\

% \linenumbers


\noindent \textbf{Acknowledgement}
\\
Fatma Yerlikaya-\"Ozkurt is supported by the TUBITAK Domestic Doctoral Scholarship Program.


%% References


% \bibliographystyle{model3a-num-names}
% \bibliography{<your-bib-database>}

%% not want to use model3a-num-names.bst.


\end{document}